\documentclass[10pt,reqno]{amsart}

\usepackage{amsmath,amssymb,amsthm}

\def\dis
{\displaystyle}

\def\R{{\mathbb R}}

\def\S{{\mathcal S}}
\def\virgp{\raise 2pt\hbox{,}}

\def\1{{\rm 1\mskip-4.5mu l} }

\def\<{\langle}
\def\>{\rangle}

\def\({\left(}
\def\){\right)}

\def\Tend#1#2{\mathop{\longrightarrow}\limits_{#1\rightarrow#2}}

\def\d{{\partial}}

\def\si{{\sigma}}
\def\om{\omega}
\def\w{{\tt w}}

\theoremstyle{plain}
\newtheorem{theo}{Theorem}[section]
\newtheorem{lem}[theo]{Lemma}

\newtheorem{prop}[theo]{Proposition}

\theoremstyle{definition}

\theoremstyle{remark}
\newtheorem{rema}[theo]{Remark}

\numberwithin{equation}{section}

\begin{document}

\title[Global existence results for NLS with quadratic
potentials]{Global existence results for nonlinear Schr\"odinger 
equations with quadratic potentials}
\author[R. Carles]{R{\'e}mi Carles}
\address{MAB, UMR CNRS 5466\\
Universit{\'e} Bordeaux 1\\ 351 cours de la Lib{\'e}ration\\ 33~405 Talence
cedex\\ France}
\email{Remi.Carles@math.u-bordeaux1.fr}
\urladdr{http://www.math.u-bordeaux1.fr/\textasciitilde carles}
\begin{abstract}
We prove that no finite time blow up can occur for nonlinear Schr\"odinger 
equations with quadratic potentials, provided that the potential has a
sufficiently strong repulsive component. This is not obvious in
general, since the energy associated to the linear equation is not
positive. The proof relies essentially on two arguments: global
in time Strichartz estimates, and a refined analysis of the linear
equation, which makes it possible to use continuity arguments and to
control the nonlinear effects. 
\end{abstract}
\maketitle

\section{Introduction}
\label{sec:intro}

We consider the nonlinear Schr\"odinger equation on $\R^n$, 
\begin{equation}\label{eq:nlsp}
i\d_t u +\frac{1}{2}\Delta u = V(x)u+ \lambda |u|^{2\si}u\quad ; \quad
u\big|_{t=0} = u_0\, ,
\end{equation}
when the potential $V$ is quadratic in $x$, or more generally when $V$
is a second order polynomial. The general assumptions we make are the
following: the space variable $x$ is in $\R^n$, $\lambda \in \R$,
$\si>0$ with $\si<\frac{2}{n-2}$ 
if $n\ge 3$ (the nonlinearity is $H^1$--sub-critical), and 
\begin{equation*}
u_0 \in \Sigma := \left\{ f \in \S'(\R^n)\ ;\ \left\| f\right\|_\Sigma
:= \|f\|_{L^2(\R^n)} + \|\nabla f\|_{L^2(\R^n)}+ \|x
f\|_{L^2(\R^n)}<+\infty \right\}\, . 
\end{equation*}
This space is very natural when one studies the case of the harmonic
potential $V(x)\equiv |x|^2$, see e.g. \cite{CazCourant}. When the
nonlinearity is $L^2$--sub-critical ($\si<2/n$), one can even consider
initial data in $L^2$ only. 
When $V\in C^\infty(\R^n;\R)$ is sub-quadratic ($\d^\alpha V\in
L^\infty(\R^n)$ for $|\alpha|\ge 2$), it is easy to prove existence
and uniqueness of solutions to \eqref{eq:nlsp} in $\Sigma$ locally in
time, without making any assumption on the sign of $V$ (see
e.g. \cite{CazCourant}). This follows for instance from the 
fact that dispersive estimates for the linear equation ($\lambda =0$)
are available for small time intervals (using perturbation arguments, one
can construct a parametrix, see \cite{Fujiwara79,Fujiwara}). One can
then apply a fixed point argument on the Duhamel's formula associated
to \eqref{eq:nlsp}, on a
closed system involving $u$, $\nabla_x u$ and $xu$. 

When $V$ is exactly the above
isotropic harmonic potential, it is easy to prove global existence
 of the solution of \eqref{eq:nlsp} in $\Sigma$, when $\lambda>0$ for
instance, thanks to the conservations of mass and energy:
\begin{equation}\label{eq:conservations}
\begin{aligned}
&\left\| u(t)\right\|_{L^2} \equiv \left\|
u_0\right\|_{L^2}\, ,\\
& E_V := \frac{1}{2}\left\| \nabla_x
u(t)\right\|_{L^2}^2 +\frac{\lambda}{\si +1}\left\|
u(t)\right\|_{L^{2\si +2}}^{2\si +2}
+\int_{\R^n}V(x)\left|u(t,x)\right|^2 dx =\text{const.}
\end{aligned}
\end{equation}
When $V$ is nonnegative (or just bounded from below), and $\lambda>0$,
these conservations yield 
\emph{a priori} estimates from which global existence follows. 

The opposite case is when $V(x)=-|x|^2$ is the repulsive harmonic
potential. The quadratic case is critical on the one hand to ensure
that the operator $-\Delta +V$ is essentially self-adjoint on
$C_0^\infty(\R^n)$: for $-\Delta -|x|^4$, this property
fails (classical trajectories can reach infinite speed, see
\cite{ReedSimon2,Dunford}).  On the other hand, it was 
proved in \cite{CaSIMA} that despite this first negative impression,
the repulsive harmonic potential tends to encourage global
existence. Of course, such results cannot be proved from the single
conservations \eqref{eq:conservations}, since the linear energy
($\lambda =0$) is not even a positive functional. Global existence in
the case $\lambda >0$ for instance stems from a conservation law,
which can be viewed as the analog of the pseudo-conformal conservation
law of the nonlinear Schr\"odinger equation with no potential (see
\cite{CaSIMA} for more details). 

The above mentioned law seems to hold only for \emph{isotropic}
repulsive harmonic potentials. If $n\ge 2$ and $V(x) = -x_1^2$, then
no global existence result seems to be available, even in the case
$\lambda>0$. The aim of this paper is to give sufficient conditions on
the potential $V$ so that the solution $u$ is global in time. This
issue is part of the more general framework to understand the interaction
between the linear dynamics generated by $-\Delta+V$, and nonlinear
effects.  

As was already exploited in different contexts \cite{CaAHP,CaSIMA,CM},
the fundamental solution for the linear problem is given explicitly
when the potential $V$ is a second order polynomial. This ``miracle''
is known as \emph{Mehler's formula} (see \cite{Feyn,HormanderQuad}),
whose expression is given below (see \eqref{eq:mehler}), and which is
closely related to the fact that for such potentials, everything is
known about classical trajectories. On the other hand, there is a  gap between
this nice framework and a more general case (see
\cite{CM} for a more quantified discussion). This is why our study is
restricted to potentials which are second order polynomials. 

After
reduction (see \cite{CM}), we may suppose that 
\begin{equation*}
V(x) =\sum_{j=1}^n \( \delta_j\frac{\om_j^2}{2}x_j^2 + b_j x_j\)\ 
, \ n\ge 1\, ,
\end{equation*}
where $\om_j>0$, $\delta_j\in \{-1,0,1\}$ and $\delta_j b_j=0$ for any
 $j$. As noticed in \cite{CN}, the Avron--Herbst formula makes it
 possible to remove the linear terms without affecting the global
 existence issue: we assume $b_j=0$ for any $j$, and $V$ is of the form
\begin{equation}\label{eq:V}
V(x) =\sum_{j=1}^n  \delta_j\frac{\om_j^2}{2}x_j^2 \ 
; \ n\ge 1,\ \om_j>0,\ \delta_j\in \{-1,0,1\}\, ,\quad \delta_1=-1\, .
\end{equation}
The last assumption means that we do not consider positive potentials,
for which several results are available (see e.g. \cite{CazCourant}). 
We denote
\begin{equation*}
H_V = -\frac{1}{2}\Delta +V\quad ;\quad  
U_V(t) = e^{-itH_V}\, .
\end{equation*}
Even though the energy $E_V$ defined in \eqref{eq:conservations} is
not a positive functional in the linear case ($\lambda =0$) when
$V$ satisfies \eqref{eq:V}, we can prove global existence in the $L^2$
sub-critical case $\si<\frac{2}{n}$, thanks to a careful analysis of
Strichartz estimates:
\begin{prop}\label{prop:general}
Let $n\ge 1$, $\lambda\in \R$, $\si<\frac{2}{n-2}$ if $n\ge 3$, and
$V$ satisfying 
\eqref{eq:V}.\\  
$1.$ Suppose $u_0\in
L^2(\R^n)$. 
\begin{itemize}
\item[(i)] If $\si<\frac{2}{n}$, then \eqref{eq:nlsp} has a unique solution 
\begin{equation*}
u\in C\cap L^\infty \(\R;L^2(\R^n)\) \cap L^q_{\rm loc}\( \R ;
L^{2\si+2}(\R^n)\)\, ,
\end{equation*}
where $q=\frac{4\si+4}{n\si}$. In addition, the $L^2$-norm of
$u(t,\cdot)$ is independent of time.
\item[(ii)] If $\si \le\frac{2}{n}$, then there
exists $\delta =\delta (\si,n,|\lambda|,V)>0$ such that
if $\|u_0\|_{L^2} <\delta$, then \eqref{eq:nlsp} has a unique global
solution $u\in C(\R;L^2)\cap L^{2+2\si}(\R\times\R^n)$.  In
addition, the $L^2$-norm of
$u(t,\cdot)$ is independent of time, and there is scattering: there 
exist unique $u_-, 
u_+\in L^2$  such that 
\begin{equation*}
\left\| U_V(-t)u(t)-u_\pm\right\|_{L^2} \Tend
t {\pm\infty} 0\, .
\end{equation*}
\end{itemize}
$2.$ Suppose $u_0\in \Sigma$.  
\begin{itemize}
\item[(i)] If $\si<\frac{2}{n}$, then the
 solution $u\in L^\infty \(\R;L^2(\R^n)\)$ of \eqref{eq:nlsp} is also
in $C(\R;\Sigma)$, and satisfies \eqref{eq:conservations}. 
\item[(ii)] There
exists $\eta =\eta (\si,n,|\lambda|,V)>0$ such that
if $\|u_0\|_\Sigma <\eta$, then \eqref{eq:nlsp} has a unique global
solution $u\in C(\R;\Sigma)$.  In addition, there is scattering: there
exist unique $u_-, 
u_+\in \Sigma$  such that 
\begin{equation*}
\left\| U_V(-t)u(t)-u_\pm\right\|_\Sigma \Tend
t {\pm\infty} 0\, .
\end{equation*}
\end{itemize}
\end{prop}
\begin{rema}\label{rem:sousquad}
When $\si<2/n$, the points 1.i and 2.i 
hold more generally when $V$ 
is a sub-quadratic potential (see Section~\ref{sec:general}). 
\end{rema}
For $\si\ge \frac{2}{n}$ and data not necessarily small, the basic example
which we treat is when $V$ is a horse shoe:  
\begin{equation}\label{eq:horse}
i\d_t u +\frac{1}{2}\Delta u = \frac{1}{2}\left(-\om^2_1
x_1^2 + \om_2^2 x_2^2\right)u+\lambda |u|^{2\si}u\quad ; \quad
u\big|_{t=0}=u_0\, . 
\end{equation}
In Section~\ref{sec:exist}, we prove the following result:
\begin{theo}\label{theo:baby}
Let $n\ge 2$, $\lambda\in \R$, $\si\ge \frac{2}{n}$ with
$\si<\frac{2}{n-2}$ if $n\ge  
3$, and $u_0\in \Sigma$. Then there exists $\Lambda=\Lambda (n,
\si,|\lambda|, \|u_0\|_\Sigma)$ such that for 
\begin{equation}\label{eq:domin}
\om_1 \ge \Lambda (1+\om_2)+ \frac{2\si^2}{2-(n-2)\si}(1+\om_2)\ln
(1+\om_2)\, , 
\end{equation}
the solution $u$ to \eqref{eq:horse} is global in time: $u\in
C(\R;\Sigma)$. Moreover, there is scattering: there exist unique $u_-,
u_+\in \Sigma$  such that 
\begin{equation*}
\left\| U_V(-t)u(t)-u_\pm\right\|_\Sigma \Tend
t {\pm\infty} 0\, .
\end{equation*}
\end{theo}
The  theorem has a heuristic explanation. The $\om_1$
term corresponds to a repulsive force: the larger $\om_1$, the
stronger the repulsive force. On the other hand, the $\om_2$ term
tends to confine the solution. There is a competition between these 
two effects. The above statement means that if the repulsive force
dominates the confining one, then the solution
cannot blow up in finite time, provided that in addition,
$\om_1\gg 1$. The latter assumption is in the same spirit as the
results of \cite{CaSIMA} in the case $\lambda<0$ and $\si\ge 2/n$. 
The second part of the theorem means that nonlinear effects become
negligible for large times. Theorem~\ref{theo:baby} yields a
\emph{sufficient} condition for global existence: it is not clear
whether the assumption \eqref{eq:domin} is sharp or not (see
Section~\ref{sec:discussion} for a discussion).  

Theorem~\ref{theo:baby} has a straightforward generalization, with a
similar proof which we omit, because it involves heavier notations and
bears no new difficulty. 

\begin{theo}\label{theo:gene}
Let $n\ge 2$, $\lambda\in \R$, $\si\ge \frac{2}{n}$ with
$\si<\frac{2}{n-2}$ if $n\ge  
3$, and $u_0\in \Sigma$. Let $V$ be of the form \eqref{eq:V}, and
denote
\begin{equation*}
\om_\pm = \max \left\{ \om_j\ ;\ \delta_j=\pm 1\right\}\quad (\om_+
=0\text{ if there is no }\delta_j=+1)\, .
\end{equation*}
Then there exists $\Lambda=\Lambda (n,
\si, |\lambda|,\|u_0\|_\Sigma)$ such that for 
\begin{equation*}
\om_- \ge \Lambda (1+\om_+)+ \frac{2\si^2}{2-(n-2)\si}(1+\om_+)\ln
(1+\om_+)\, , 
\end{equation*}
the solution $u$ to \eqref{eq:nlsp} is global in time, $u\in
C(\R;\Sigma)$, and there is scattering.
\end{theo}
As above, the statement can be summarized as follows: if the repulsive
force is sufficiently strong compared to other effects (linear
confinement is overcome if $\om_-\gg \om_+$, nonlinear effects are 
overcome if $\om_-\gg 1$), then the solution is global and the
nonlinearity can be viewed as a (short range) perturbation for large
times. 

The paper is organized as follows. In Section~\ref{sec:strichartz}, we
recall Strichartz estimates, and notice that they hold globally in
time for $U_V$, when $V$ is of the form \eqref{eq:V}. This
follows from a simple remark 
after the proof of non-endpoint estimates in \cite{KT}. We deduce
Proposition~\ref{prop:general} in Section~\ref{sec:general}, and
Theorem~\ref{theo:baby} in Section~\ref{sec:exist}. In
Section~\ref{sec:discussion}, we finally discuss the above results
further into details.

\section{Strichartz estimates}
\label{sec:strichartz}

As recalled in the introduction, Strichartz estimates are the modern
tool to study (among others) Schr\"odinger equations. More precisely,
assume that $U$ is a \emph{$L^2$-unitary group}\footnote{With Schr\"odinger
equations with time-independent potentials in mind, this is natural.}
such that: 
\begin{equation}\label{eq:dispKT}
\left\|U(t)\right\|_{L^1\to 
L^\infty}  
\lesssim |t|^{-\gamma}\, .
\end{equation}
In the case of Schr\"odinger equations, we have $\gamma =\frac{n}{2}$,
but it does not cost much to consider a general $\gamma$. Following
\cite{KT}, we say that a pair
$(q,r)$ is \emph{sharp $\gamma$-admissible} if $q,r\ge 2$,
$(q,r,\gamma)\not =(2,\infty,1)$ and
\begin{equation*}
\frac{1}{q}+\frac{\gamma}{r}=\frac{\gamma}{2}\, \cdot
\end{equation*}
The following result is proved in \cite{KT} (see references therein
for earlier proofs). For any sharp $\gamma$-admissible pairs $(q,r)$
and $(\widetilde q,\widetilde r)$, there exist $C_q$ and
$C_{q,\widetilde q}$ such that for any time interval $I\ni 0$, 
\begin{align}
\left\| U(t)f\right\|_{L^q(I;L^r(\R^n))} &\le C_q \left\|
f\right\|_{L^2(\R^n)}\, ,\label{eq:strihom}\\
\left\| \int_0^t U(t-s)F(s)ds \right\|_{L^q(I;L^r(\R^n))} &\le
C_{q,\widetilde q} \left\| F\right\|_{L^{\widetilde
q'}(I;L^{\widetilde 
r'}(\R^n))}\, , \label{eq:striinhom}
\end{align}
where $a'$ stands for the H\"older conjugate exponent of $a$. 

In order to
explain the existence of this section, recall Mehler's formula
(\cite{Feyn,HormanderQuad}). If $V$ is of the form \eqref{eq:V}, then
denoting $H_V=-\frac{1}{2}\Delta +V$, we have:
\begin{equation}\label{eq:mehler}
U_V(t)f:= e^{-it H_V}f=
\prod_{j=1}^n
\left(\frac{1}{2i\pi g_j(t)}\right)^{1/2}
\int_{\R^n}e^{iS(t,x,y)}f(y)dy\, ,
\end{equation}
where
\begin{equation*}
S(t,x,y)= \sum_{j=1}^n \frac{1}{g_j(t)}\(
\frac{x_j^2 +y_j^2}{2}h_j(t) -x_j
y_j\),
\end{equation*}
and the functions $g_j$ and $h_j$, related to the classical
trajectories, are given by:
\begin{equation}\label{eq:rays} 
\begin{aligned}
\(g_j(t),h_j(t)\) =&\left\{
\begin{aligned}
\( \frac{\sinh (\om_j t)}{\om_j},\cosh (\om_j t)\)\, ,\ &\text{ if
}\delta_j=-1\, ,\\ 
\(t,1\)\ ,\ &\text{ if }\delta_j=0 \, ,\\  
\(\frac{\sin (\om_j t)}{\om_j}, \cos(\om_j t)\)\, ,\ &\text{ if }
\delta_j=+1\, .
\end{aligned}
\right.
\end{aligned}
\end{equation}
Recall that if there exists $\delta_j=+1$, then $e^{-it{H}_V}$ has
some singularities, periodically in time  
(see e.g. \cite{KR}). This affects the 
above formula with phase factors we did not write (which can be
incorporated in the definition of $(ig_j(t))^{1/2}$). However, these
singularities may prevent the existence of \emph{global in time}
Strichartz estimates. Indeed, we obviously have 
\begin{equation*}
\left\|U_V(t)\right\|_{L^1\to 
L^\infty}  
\lesssim |t|^{-n/2} \quad \text{for }|t|\le \delta \, .
\end{equation*}
This yields local in time Strichartz estimates: the above constants
$C_q$ and $C_{q,\widetilde q}$ depend on the time interval $I$, and  
they may blow up on unbounded time intervals. This is obvious in the
case of the isotropic harmonic potential, since eigenfunctions yield
non-dispersive solutions (this can also be read from Mehler's
formula). On the other hand, the exponential decay provided by a
repulsive component of the potential (that is, there exists at least
one $\delta_j=-1$) suggests that the singularities of confining forces
may be balanced. We show that this is the case. 

To see this, we simply remark that the above Strichartz estimates are
still valid if we replace the assumption \eqref{eq:dispKT} by 
\begin{equation}\label{eq:dispgen}
\left\|U(t)\right\|_{L^1\to 
L^\infty}  
\le \w(t)^{\gamma}\, ,\quad \text{with }\w\ge 0\ \text{and }\w
\in L^1_w(\R)\, ,
\end{equation}
the weak $L^1$ space\footnote{The assumption \eqref{eq:dispKT}
is $\w(t)=c|t|^{-1}$ which is the basic example for a function in
$L^1_w(\R)$.}. This is straightforward, since the proof in 
\cite{KT} is actually valid with this relaxed assumption, for non-endpoint
estimates. For the convenience of the reader, we recall the argument,
with $I=\R$.

First, by duality, \eqref{eq:strihom} is equivalent to:
\begin{equation}\label{eq:6}
\left\| \int U(-s)F(s)ds\right\|_{L^2(\R^n)}\lesssim \left\|
F\right\|_{L^{q'}(\R;L^{r'}(\R^n))}\, . 
\end{equation}
By the $TT^*$ method, this is equivalent to 
\begin{equation*}
\left| \iint \<U(-s)F(s), U(-t)G(t)\>dsdt\right|\lesssim \left\|
F\right\|_{L^{q'}(\R;L^{r'}(\R^n))}\left\|
G\right\|_{L^{q'}(\R;L^{r'}(\R^n))}\, . 
\end{equation*}
By symmetry, it suffices to prove 
\begin{equation}\label{eq:14}
|T(F,G)|\lesssim \left\|
F\right\|_{L^{q'}_t L^{r'}_x}\left\|
G\right\|_{L^{q'}_t L^{r'}_x}\, , 
\end{equation}
where $\dis T(F,G) = \iint_{s<t}
\<U(-s)F(s), U(-t)G(t)\>dsdt$.
Since $U$ is unitary on $L^2$,
\begin{equation*}
\left|  \<U(-s)F(s), U(-t)G(t)\>\right|\le \left\|
F(s)\right\|_{L^2_x} \left\|G(s)\right\|_{L^2_x} \, .
\end{equation*}
We infer from the dispersive estimate \eqref{eq:dispgen} that 
\begin{align*}
\left|  \<U(-s)F(s), U(-t)G(t)\>\right|&= \left|  \<U(t-s)F(s),
G(t)\>\right|\\
&\le  \w(t-s)^\gamma 
\left\|F(s)\right\|_{L^1_x} \left\|G(s)\right\|_{L^1_x} \, .
\end{align*}
By interpolation,
\begin{equation*}
\left|  \<U(-s)F(s), U(-t)G(t)\>\right|\le \w(t-s)^{2/q}\left\|
F(s)\right\|_{L^{r'}_x} \left\|G(s)\right\|_{L^{r'}_x} \, ,
\end{equation*}
since $(q,r)$ is sharp $\gamma$-admissible (and is not an
endpoint). Then \eqref{eq:14} 
follows from Hardy--Littlewood--Sobolev inequality (see
e.g. \cite{Stein70} or \cite[Sect.~4.3]{LiebLoss}). This 
proves the homogeneous estimate \eqref{eq:strihom}. 

For the inhomogeneous case, consider two sharp $\gamma$-admissible
pairs $(q,r)$ 
and $(\widetilde q,\widetilde r)$.
By duality, \eqref{eq:striinhom} is equivalent to:
\begin{equation}\label{eq:32}
\left| T(F,G)\right| \lesssim \left\| F\right\|_{L^{q'}_tL^{r'}_x}
\left\| G\right\|_{L^{\widetilde q'}_tL^{\widetilde r'}_x} \, .
\end{equation}
We have
\begin{equation*}
\left| T(F,G)\right| \le \(\sup_{t\in\R}\left\|\int_{s<t}U(-s)F(s)ds
\right\|_{L^2_x}\) 
\left\| G\right\|_{L^1_tL^2_x}\, ,
\end{equation*}
and when $(\widetilde q,\widetilde r)=(\infty,2)$, \eqref{eq:32}
follows from \eqref{eq:6}. Similarly, \eqref{eq:32} holds when
$( q,r)=(\infty,2)$. From \eqref{eq:14}, one has \eqref{eq:32} when
$(q,r)=(\widetilde q,\widetilde r)$. The general case
\eqref{eq:striinhom} follows by interpolation between these three
cases. We have precisely:
\begin{lem}\label{lem:strichartz}
Let $(U(t))_{t\in\R}$ be a unitary group on $L^2(\R^n)$, satisfying
the dispersive estimate \eqref{eq:dispgen}. Then for any $T\in
\overline{\R}_+$, and any sharp 
$\gamma$-admissible pairs $(q,r)$ and $(\widetilde q,\widetilde r)$,
we have
\begin{align*}
\left\| U(t)f\right\|_{L^q(]-T,T[;L^r)} &\le c_q
\left\|\w\1_{]-2T,2T[}\right\|_{L^1_w}^{1/q} \left\| 
f\right\|_{L^2}\, ,\\
\left\| \int_0^t U(t-s)F(s)ds \right\|_{L^q(]-T,T[;L^r)} &\le
C_{q,\widetilde q}
\left\|\w\1_{]-2T,2T[}\right\|_{L^1_w}^{1/q +1/\widetilde q}
\left\| F\right\|_{L^{\widetilde q'}(]-T,T[;L^{\widetilde
r'})}\, .
\end{align*}
\end{lem}

From now on, we shall simply write that a pair $(q,r)$ is
admissible when it is sharp $\frac{n}{2}$-admissible.

In the case of a quadratic potential \eqref{eq:V}, this yields global
in time Strichartz estimates, since $\delta_1=-1$. From Mehler's
formula \eqref{eq:mehler}, we have, for some $\delta>0$,
\begin{equation*}
\left\|U_V(t)\right\|_{L^1\to 
L^\infty}  
\lesssim |t|^{-n/2} \quad \text{for }|t|\le \delta \, .
\end{equation*}
Now for $|t|>\delta$, the ``worst'' possible case is when, say,
$\delta_1=-1$ and $\delta_j =+1$ for $j\ge 2$. Then
\begin{equation*}
\left\|U_V(t)\right\|_{L^1\to 
L^\infty}  
\lesssim \(e^{-\om_1 |t|}\prod_{j=2}^n \frac{1}{|\sin (\om_j t)|}\)^{1/2}
\quad \text{for }|t|> \delta \, . 
\end{equation*}
Therefore, \eqref{eq:dispgen} is satisfied with 
\begin{equation*}
\w(t) = \text{const.}\(\frac{1}{|t|}\1_{|t|\le \delta} + \(e^{-\om_1
|t|}\prod_{j=2}^n 
\frac{1}{|\sin (\om_j t)|}\)^{1/n} \1_{|t|> \delta}\)\,.  
\end{equation*}
The first part is obviously in $L^1_w(\R)$. The second part is in
$L^1(\R)$; this follows from H\"older's inequality, and the fact that 
\begin{equation*}
t\mapsto \(\frac{e^{-|t|}}{|\sin t|}\)^{1/n} \in L^{n-1}(\R)\, .
\end{equation*}
Therefore, $\w \in L^1_w(\R)$, and we have global in time Strichartz
inequalities. In the case of \eqref{eq:horse}, that is when
$V(x)=\frac{1}{2}(-\om_1^2x_1^2 +\om_2^2x_2^2)$, we have:
\begin{align}
&\left\|\w\1_{]-\frac{\pi}{2\om_2},\frac{\pi}{2\om_2}[}\right\|_{L^1_w} 
\lesssim 1\quad \text{(always)}\, ,\label{eq:wpetit} \\
&\|\w\|_{L^1_w} \lesssim 1 + \( \frac{\om_1}{\om_2}\)^{\frac{1}{n(n-1)}}
e^{-\frac{\om_1}{\om_2} \frac{\pi}{2n(n-1)}}\, ,\text{ if
}\om_1\ge \om_2\, .\label{eq:wpetit2}
\end{align}
\eqref{eq:wpetit} is straightforward; \eqref{eq:wpetit2} follows from
H\"older's inequality, writing $\w$ as the product of $n-1$ terms, for
$|t|>\frac{\pi}{2\om_2}$, 
\begin{equation*}
\w(t)  = \frac{1}{\sqrt{2\pi}} \( \(\frac{\om_1}{|\sinh (\om_1
t)|}\)^{\frac{1}{n-1}}\frac{\om_2}{|\sin (\om_2 t)|}\) \times 
\prod_{j=1}^{n-2} \(\frac{1}{|t|^{1/n}}\(\frac{\om_1}{|\sinh (\om_1
t)|}\)^{\frac{1}{n-1}}\)\, .  
\end{equation*}

\begin{rema}\label{rem:dimension}
For $d>0$, define the dispersive rate $\w_d$ by 
$\w(t)^{n/2}=\w_d(t)^{d/2}$.
It is easy to check that $\w_d\in L^1_w(\R)$ provided that $d\ge n$
(and  $\w_d\in L^1 (|x|>1)$). Therefore,
Lemma~\ref{lem:strichartz} shows that Strichartz inequalities hold
with the same admissible pairs as in ``space dimension $d$'' (even
though $d$ needs not be an integer). 
\end{rema}

\begin{rema}
Notice that we have global in time Strichartz estimates in cases where
there exist trapped trajectories. Consider for instance 
\eqref{eq:horse} with $n=2$ and $\om_1=\om_2=1$. In general, classical
trajectories solve $\ddot x +\nabla V(x)=0$; this yields in the
present 
case:
\begin{align*}
x_1(t) &= x_1(0)\cosh t +\xi_1(0)\sinh
 t= 
\frac{e^{ t}}{2}\( x_1(0)+\xi_1(0)\)
+\frac{e^{- t}}{2}\( x_1(0)-\xi_1(0)\) \, ,\\
x_2(t)& = x_2(0)\cos t +\xi_2(0)\sin t\, . 
\end{align*}
If $x_1(0)+\xi_1(0)=0$, then the trajectory is trapped
in the future, but not in the past in general. If we
suppose in addition that $x_1(0)=\xi_1(0)=0$, then we have a
trajectory which is trapped in the past \emph{and} in the future (and
nontrivial if $\xi_2(0)\not =0$). 
Compare with the 
results of \cite{CKS95,Doi00}; 
 it is proved that smoothing effects occur in the future provided
that the classical trajectories are not trapped in the past. However,
the results of \cite{CKS95} include potentials which grow at most
linearly in $x$, and \cite{Doi00} does not consider the case of
potentials. On the other hand,
smoothing effects yield another approach to prove Strichartz estimates
(see e.g. \cite{ST02,BGTIHP,BP}). In our
case, there exist trajectories trapped in the past and in the
future, but global in time Strichartz estimates are available. 
It seems
that the link between classical trajectories and (global in time)
Strichartz estimates remains to be clarified. 
\end{rema}

\section{Proof of Proposition~\ref{prop:general}}
\label{sec:general}

In the previous section, we saw that when $V$ is of the form
\eqref{eq:V}, $U_V$ satisfies global in time Strichartz estimates. The
constants in these inequalities may depend on $\om_1,\ldots,\om_n$,
but this is not important in view of
Proposition~\ref{prop:general}. As a matter of fact, global in time
Strichartz estimates are really needed only for the third point of the
proposition.

The first part of Proposition~\ref{prop:general} is
straightforward: 
one can mimic the proof given in the case of the nonlinear
Schr\"odinger equation with no potential ($V\equiv 0$ in
\eqref{eq:nlsp}, see
\cite{TsutsumiL2,CW89,CazCourant}). We recall the main argument. Duhamel's
formula writes  
\begin{equation}\label{eq:duhamel}
u(t)=U_V(t)u_0 -i\lambda \int_0^t U_V(t-s)\(
|u|^{2\si}u \)(s)ds.
\end{equation}
Define $F(u)(t)$ as the right hand side of \eqref{eq:duhamel}. The
idea is to use a fixed point argument in the space given in
Proposition~\ref{prop:general}. Introduce the following Lebesgue
exponents:  
\begin{equation}\label{eq:exponents}
r=2\si +2\quad ;\quad
q=\frac{4\si+4}{n\si}\quad ;\quad k=
\frac{2\si(2\si+2)}{2-(n-2)\si}\,\cdot
\end{equation}
Then $(q,r)$ is the admissible pair of the proposition, and 
\begin{equation*}
\frac{1}{r'}=\frac{2\si}{r}+\frac{1}{r}\quad ; \quad
\frac{1}{q'}=\frac{2\si}{k}+\frac{1}{q}\, \cdot
\end{equation*} 
The main remark to prove the first point of
Proposition~\ref{prop:general} is that if $\si<\frac{2}{n}$, we
have $\frac{1}{q}<\frac{1}{k}$, and H\"older's inequality in time
yields
\begin{equation*}
\| u\|_{L^k(I;L^r)} \le |I|^{1/k-1/q} \| u\|_{L^q (I;L^r)} =
|I|^{\frac{(2-n\si)(\si+1)}{2\si(2\si+2)}} \| u\|_{L^q 
(I;L^r)}\, .
\end{equation*}
The positive power of $|I|$ yields contraction in $L^\infty L^2\cap
L^qL^r$ for small time intervals, and the conservation of the $L^2$
norm of the solution shows global existence at the $L^2$ level. This
proves 1.i.

If $\si=\frac{2}{n}$, then $k=q$, and Lemma~\ref{lem:strichartz} yields
\begin{equation*}
\| u\|_{L^{2+\frac{4}{n}}(I\times\R^n)}\le C\|u_0\|_{L^2} + C \|
u\|_{L^{2+\frac{4}{n}}(I\times\R^n)}^{1+\frac{4}{n}}\, ,
\end{equation*}
for some constant $C$ \emph{independent of the time interval} $I$. The idea
is then to use a bootstrap argument, for $\|u_0\|_{L^2}$ sufficiently
small (see \cite{CW89,CazCourant} for details). When
$\si<\frac{2}{n}$, one can use the same approach, thanks to
Remark~\ref{rem:dimension}. Define $d= \frac{2}{\si}>n$. Then the
nonlinearity is $L^2$--critical in ``space dimension $d$'', and the
proof in \cite{CW89} can be applied. This
completes the proof of the first part of
Proposition~\ref{prop:general}. Note that in the small data case, we
used the fact that we have global in time Strichartz estimates, due to
the repulsive character of the potential, $\delta_1=-1$ in
\eqref{eq:V}. 

To prove the second point of Proposition~\ref{prop:general}, we
restrict to the case of \eqref{eq:horse}. This gives all the arguments
for the general case, and prepares the proof of
Theorem~\ref{theo:baby}. 

Introduce the operators
\begin{align*}
J_1(t) = \om_1 x_1 \sinh(\om_1 t) +i\cosh (\om_1 t)\d_1\quad &; \quad
H_1(t) = x_1 \cosh (\om_1 t) +i\frac{\sinh(\om_1 t)}{\om_1}\d_1 \, ,\\
J_2(t) = -\om_2 x_2 \sin(\om_2 t) +i\cos (\om_2 t)\d_2\quad &; \quad
H_2(t) = x_2 \cos (\om_2 t) +i\frac{\sin(\om_2 t)}{\om_2}\d_2 \, .
\end{align*}
We define $J = (J_k)_{1\le k\le n}$ and $H = (H_k)_{1\le k\le n}$,
where if $n\ge 3$, 
\begin{equation*}
J_k(t) = i\d_k \quad ; \quad H_k(t) = x_k +it\d_k\quad \text{for }k\ge
3\, .
\end{equation*}
We have the weighted Gagliardo--Nirenberg inequalities:
\begin{lem}\label{lem:GN}
Let $2\le p< \frac{2n}{n-2}$. There exists $C_p$
independent of $\om_1,\om_2>0$ such that for any $f\in \Sigma$, 
\begin{equation*}
\left\| f\right\|_{L^p} \leq C_p\(\frac{\left\|
J_1(t)f\right\|_{L^2}}{\cosh(\om_1
t)}\)^{\frac{\delta(p)}{n}}\( \frac{\left\|
J_2(t)f\right\|_{L^2}}{|\cos(\om_2
t)|}\)^{\frac{\delta(p)}{n}} \left\|
f\right\|_{L^2}^{1-\delta(p)}\prod_{j=3}^n \left\| 
\d_j f\right\|_{L^2}^{\frac{\delta(p)}{n}}   ,
\end{equation*}
where $\delta(p)= n\(\frac{1}{2}-\frac{1}{p}\)$. We also have
\begin{equation*}
\left\| f\right\|_{L^p} \leq C_p\(\frac{\left\|
J_1(t)f\right\|_{L^2}}{\cosh(\om_1
t)}\)^{\frac{\delta(p)}{n}}\( \om_2 \frac{\left\|
H_2(t)f\right\|_{L^2}}{|\sin(\om_2
t)|}\)^{\frac{\delta(p)}{n}} \left\|
f\right\|_{L^2}^{1-\delta(p)}\prod_{j=3}^n \left\| 
\d_j f\right\|_{L^2}^{\frac{\delta(p)}{n}}   ,
\end{equation*}
and therefore
\begin{equation*}
\left\| f\right\|_{L^p} \lesssim\left\|
f\right\|_{L^2}^{1-\delta(p)}\(\frac{\left\|
J_1(t)f\right\|_{L^2}}{\cosh(\om_1
t)}\( \left\|
J_2(t)f\right\|_{L^2} + \om_2 \left\|
H_2(t)f\right\|_{L^2}\) \prod_{j=3}^n \left\| 
\d_j f\right\|_{L^2}\)^{\frac{\delta(p)}{n}}   .
\end{equation*}
\end{lem}
This lemma follows from the usual Gagliardo--Nirenberg inequalities
and:
\begin{equation}\label{eq:factor}
\begin{aligned}
J_1(t) &= i \cosh (\om_1 t) e^{i\om_1 \frac{x_1^2}{2}
\tanh(\om_1 t)}\d_1 \( e^{-i\om_1 \frac{x_1^2}{2}
\tanh(\om_1 t)}\ \cdot \)\, ,\\
J_2(t) &= i \cos (\om_2 t) e^{-i\om_2 \frac{x_2^2}{2}
\tan(\om_2 t)}\d_2 \( e^{i\om_2 \frac{x_2^2}{2}
\tan(\om_2 t)}\ \cdot \)\, ,\\
H_2(t) &= i \frac{\sin (\om_2 t)}{\om_2} e^{i\om_2 \frac{x_2^2}{2}
\cot(\om_2 t)}\d_2 \( e^{-i\om_2 \frac{x_2^2}{2}
\cot(\om_2 t)}\ \cdot \)\, .
\end{aligned}
\end{equation}
Note that the Galilean operators $H_k$, $k\ge 3$, share the same
property:
\begin{equation}\label{eq:factorbis}
H_k(t) = i t\, e^{i\frac{x_k^2}{2t}}\d_k \( e^{-i\frac{x_k^2}{2t}}\
\cdot \)\, ,\quad \text{for }k\ge 3\, .
\end{equation}
To prove global existence in $\Sigma$, it is sufficient to
prove that $A(t)u\in L^\infty_{\rm loc}(\R;L^2)$ for any $A\in
\{Id,J,H\}$. This follows from the formula
\begin{equation*}
\( 
\begin{array}{c}
J_j(t)\\
H_j(t)
\end{array}
\)= \( 
\begin{array}{cc}
-\delta_j\om_j^2 g_j(t) & h_j(t)\\
h_j(t)& g_j(t)
\end{array}
\)
\( 
\begin{array}{c}
x_j\\
i\d_j
\end{array}
\)\, ,\quad \forall j\ge 1\, .
\end{equation*}
The operators we use share the same
properties as those 
which are used for $\nabla_x$ in the case with no potential:
they commute with the linear part of the equation (including the
potential $V$, since they are Heisenberg derivatives -- see \cite{CM}
and the discussion therein); they yield Gagliardo--Nirenberg
inequalities; they act like derivatives on 
nonlinearities of the form $G(|z|^2)z$, from \eqref{eq:factor} and
\eqref{eq:factorbis}. 

Fix $A\in \{Id,J,H\}$. From the above arguments and Strichartz
estimates, 
\begin{align*}
\left\| A(t) F(u)\right\|_{L^q(I;L^r)} &\lesssim \left\| A(t)
\(|u|^{2\si}u\)\right\|_{L^{q'}(I;L^{r'})} \lesssim \left\| 
|u|^{2\si}A(t)u\right\|_{L^{q'}(I;L^{r'})}\\
&\lesssim \|u\|_{L^k(I;L^r)}^{2\si} \left\| 
A(t)u\right\|_{L^{q}(I;L^{r})}\, . 
\end{align*}
Since $u\in L^q_{\rm loc}(\R;L^r)$ and $\frac{1}{q}<\frac{1}{k}$, we
have $Au \in L^q_{\rm loc}(\R;L^r)$. Using Strichartz inequalities
once more, we infer that $Au \in L^\infty_{\rm loc}(\R;L^2)$.

In the small data case, one can also resume the usual proof for global
existence in $\Sigma$. Scattering for any $\si>0$ (in this small data
case) follows the same way
as global in time Strichartz estimates, since Lemma~\ref{lem:GN}
yields an exponential decay for $u$. This 
completes the proof of Proposition~\ref{prop:general}. 
\begin{rema}
As pointed out in Remark~\ref{rem:sousquad},
some results for $\si<2/n$ hold for any sub-quadratic potential $V$. If
$V\in C^\infty(\R^n;\R)$ is such that $\d^\alpha V\in L^\infty(\R^n)$
for any $|\alpha|\ge 2$, then the results of
\cite{Fujiwara79,Fujiwara} prove that the group $U_V$ has the same
dispersion has the free group $U_0$ for \emph{small times}. One can then
mimic the proof of \cite{TsutsumiL2} to infer global existence at the
$L^2$ level. Global existence in $\Sigma$ follows, considering the
closed system of estimates for $\nabla_x u$ and $xu$. On the other
hand, no scattering result (even for small data) must be expected under
these general assumptions; the case of the harmonic potential yields a
counterexample (see \cite{CaIHP} for a high-frequency analysis). 
\end{rema}

\section{Proof of Theorem~\ref{theo:baby}}
\label{sec:exist}
First, in the
same spirit as in \cite{CaSIMA}, we analyze
the local existence result, to bound from below the local existence
time, in term of the parameters $\om_1$ and $\om_2$. Then, we notice
that we obtain a time 
at which $u$ is defined and small, and for which therefore
the nonlinearity is not too strong. We consider the solution of
the linear equation (\eqref{eq:horse} with $\lambda=0$) that coincides
with $u$ at that time. A continuity argument shows that $u$ cannot move
away too much from this linear solution. Since the linear solution is
global, so is $u$. The scattering is then an easy by-product. 

The first step is:
\begin{prop}\label{prop:local}
Let $u_0\in \Sigma$. There exists $\delta=\delta(n,\si,|\lambda|,
\|u_0\|_\Sigma)>0$ independent of 
$\om_1,\om_2\ge 0$, and $\frac{\delta}{1+\om_2}\le t_0<
\frac{\pi}{8\om_2}$ such that \eqref{eq:horse} has a unique solution 
\begin{equation*}
u\in C\(]-2t_0,2t_0[;\Sigma\) \cap L^q\( ]-2t_0,2t_0[;
W^{1,2\si+2}(\R^n)\)\, ,
\end{equation*}
with $q$ given in \eqref{eq:exponents}. In addition, the conservations
\eqref{eq:conservations} hold, and there exists $C_0$ independent of
$\om_1,\om_2\ge 0$ such that
\begin{equation*}
\sup_{|t|\leq t_0}\left\|A(t)u\right\|_{L^2}\leq C_0\, ,\quad
\forall A\in\{ Id,J,H\}\, .
\end{equation*}
\end{prop}
\begin{rema}
A similar result could be proved with slightly more general data than
$u_0\in\Sigma$, following exactly the same lines. Namely, we could assume 
\begin{equation*}
u_0\in X_{1,2}:=\left\{ f\in H^1(\R^n)\ ;\ x_1f, x_2f\in
L^2(\R^n)\right\}\, . 
\end{equation*}
\end{rema}
\begin{proof}
For $|t|<\pi/(4\om_2)$,
Lemma~\ref{lem:strichartz} and \eqref{eq:wpetit} yield Strichartz
estimates with constants independent of $\om_1$ and $\om_2$, 
and the weights in the first inequality of Lemma~\ref{lem:GN} are
bounded uniformly in  $\om_1$ and $\om_2$.

We can then mimic the
usual proof (see e.g. \cite{CazCourant}), that consists 
in applying a fixed point theorem on \eqref{eq:duhamel}, using
similar arguments as in the previous section. We give the main lines
of the proof, so that the bound from below on 
$t_0$ is clear. 
Let $R:=\|u_0\|_\Sigma$. With the definition \eqref{eq:exponents},
denote 
\begin{equation*}
Y_{r}(I) := \left\{ u\in C(I;\Sigma); \ A(t)u \in 
L^q(I;L^r)\cap L^\infty(I;L^2),\ \forall A\in 
\{ Id,J,H\}\right\}\, .
\end{equation*}
We first prove that there exists $0<T<\frac{\pi}{4\om_2}$ such
that the set 
\begin{equation*}
\begin{aligned}
 X_T := \big\{ & u \in Y_r(]-T, T[); \   \|A(t)u\|_{L^2}\leq
2R\ , \ \ \forall |t|<T,\ A\in \{ Id,J,H\},\\ 
& \|A(t)u\|_{L^q(]-T, T[;L^r)}\leq 2c_q\| \w
\1_{]-\frac{\pi}{2\om_2},\frac{\pi}{2\om_2}[} \|_{L^1_w}^{1/q} 
R, \ 
\forall A\in \{ Id,J,H\} \big\} 
\end{aligned}
\end{equation*}
is stable under the map $F$, defined in the previous section (right
hand side in Duhamel's formula). The constant $c_q$ is the one that
appears in Lemma~\ref{lem:strichartz}.
We then prove that up to choosing $T$
even smaller, $F$ is a contraction on $X_T$. The natural norm on $Y_r$
is 
$$\|u\|_{Y_r}:= \sum_{A \in \{Id,J,H\}}\left(\|A(t)u\|_{L^\infty
(]-T, T[;L^2)} + 
\|A(t)u\|_{L^q(]-T, T[;L^r)} \right).$$
For any pair $(a,b)$, we
use the notation $\| f\|_{L^a_T(L^b)}= \|f\|_{L^a(]-T, T[;L^b)}$.  
Let $u \in X_T$, and $A \in \{Id,J,H\}$. We already noticed that 
\begin{equation*}
  \|A(t)F(u)\|_{L_T^\infty(L^2)}
\le R + \widetilde C  \left\|u
  \right\|^{2\si}_{L^k_T(L^r)}\left\|A(t)u
  \right\|_{L^{q}_T(L^{r})}\, ,
\end{equation*}
where $\widetilde C$ does not depend on $\om_1$ or $\om_2$. 
From Lemma~\ref{lem:GN}, we have for $T<\frac{\pi}{4\om_2}$, 
$ \left\|u
  \right\|_{L^k_T(L^r)} \lesssim R T^{1/k}$.
It follows, 
\begin{equation}
  \label{eq:op1}
  \|A(t)F(u)\|_{L^\infty_T(L^2)} \leq R + C R^{2\si+1}T^{2\si/k}\ .
\end{equation}
Use Lemma~\ref{lem:strichartz} to obtain
\begin{equation*}
\begin{split}
  \|A(t)F(u)\|_{L^q_T(L^r)}&
\le  c_q \left\| \w
\1_{]-\frac{\pi}{2\om_2},\frac{\pi}{2\om_2}[} \right\|_{L^1_w}^{1/q} R
+ C  \left\|u 
  \right\|^{2\si}_{L^k_T(L^r)}\left\|A(t)u
  \right\|_{L^{q}_T(L^{r})}\\
& \le c_q \left\| \w
\1_{]-\frac{\pi}{2\om_2},\frac{\pi}{2\om_2}[} \right\|_{L^1_w}^{1/q}R +C
  R^{2\si+1}T^{2\si/k}\ . 
\end{split}
\end{equation*}
It is now clear that if $T$ is sufficiently small, then $X_T$ is
stable under $F$.

To complete the proof of the proposition, it is enough to prove
contraction for small $T$ in 
the weaker metric $L^q(]-T,T[;L^r)$. 
We have 
\begin{align*}
    \big\| F(u_2)-F(u_1)
    \big\|_{L^q_T(L^r)}  & \leq
C\left\| \left(|u_2|^{2\si} u_2 - |u_1|^{2\si}
    u_1\right) \right\|_{L^{q'}_T(L^{r'})}\\
&\leq C\left(
    \|u_1\|^{2\si}_{L^k_T(L^r)}+
\|u_2\|^{2\si}_{L^k_T(L^r)}\right) 
\|u_2-u_1\|_{L^q_T(L^r)}.
\end{align*}
As above, we have the estimate
$\|u_j\|^{2\si}_{L^k_T(L^r)}\le C R^{2\si}T^{2\si/k}$, $j=1,2$.
Therefore,  contraction follows for $T$ sufficiently small. The only
requirements we make are $T<\frac{\pi}{4\om_2}$ and $T\le
\eta=\eta(n,\si,|\lambda|,\|u_0\|_\Sigma)$. Therefore, we can find
$\delta>0$ independent of $\om_1$ and 
$\om_2$ such that $T\ge \frac{2\delta}{1+\om_2}$, and the proposition
follows with $T=2t_0$.
\end{proof}
\begin{rema}
As pointed out by the referee, the approach of
Proposition~\ref{prop:local} could be pursued in order to complete the
proof of Theorem~\ref{theo:baby} in just one step. Indeed, it is
enough to prove that 
$\|u\|_{L^k(\R;L^r)}^{2\si}$ can be made small under the assumptions
of Theorem~\ref{theo:baby}. We saw that choosing $t_0$ sufficiently
small,
$\|u\|_{L^k_{t_0}L^r}^{2\si}$ is small. Then using
Lemma~\ref{lem:GN}, the exponential decay shows that
$\|u\|_{L^k(\{|t|>t_0\};L^r)}^{2\si}$ can be made small for $\om_1$
large enough (see the proof below). However, we chose to keep our
initial presentation, for 
we believe it gives more information on the interaction between linear
(due to the potential) and nonlinear effects. 
\end{rema}

Now the idea is that at time $t=t_0$, the $L^p$ norms of $u$ are small
for $\om_1\gg 1+\om_2$ ($p\not =2$), from Lemma~\ref{lem:GN}. In that
case, the nonlinearity becomes negligible. 
Define the ``approximate'' solution $v$ by
\begin{equation}\label{eq:v}
i\d_t v +\frac{1}{2}\Delta v = \frac{1}{2}\left(-\om^2_1
x_1^2 + \om_2^2 x_2^2\right)v\quad ; \quad
v\big|_{t=t_0}=u\big|_{t=t_0}\, ,
\end{equation}
where $t_0$ stems from Proposition~\ref{prop:local}. We also define the
error $w=u-v$. Since $v$ is defined globally, $u$ exists globally in
time if and only if $w$ does. 

We prove Theorem~\ref{theo:baby} for positive times. The proof for
negative times is similar. 
The remainder $w$ solves
\begin{equation}\label{eq:w}
i\d_t w +\frac{1}{2}\Delta w = \frac{1}{2}\left(-\om^2_1
x_1^2 + \om_2^2 x_2^2\right)w + \lambda |u|^{2\si}u\quad ; \quad
w\big|_{t=t_0}=0\, .
\end{equation}
We note that since the operators $J$ and $H$ commute with the linear
part of \eqref{eq:horse}, the conservation of mass for Schr\"odinger
equations yields a constant $C_1$ independent of
$\om_1,\om_2>0$ such that 
\begin{equation}\label{eq:estimv}
\sum_{A\in\{Id, J,K\} }
\left\|A(t)v\right\|_{L^2} \leq C_1\, , 
\end{equation}
uniformly in $t\in\R$. In view of Theorem~\ref{theo:baby}, we may
assume $\om_1\ge \om_2$, so 
that Lemma~\ref{lem:strichartz} and \eqref{eq:wpetit2} yield global in
time Strichartz estimates, independent of $\om_1$ and $\om_2$. 
From Proposition~\ref{prop:local}, there exists $t_1>0$ such that $w$
satisfies the same estimates as $v$ on $[t_0,t_0+t_1]$. So long as this
holds, Strichartz estimates, together with H\"older's
inequality, yield:
\begin{equation}\label{eq:reste} 
\begin{aligned}
\|w\|_{L^q(I;L^r)} &\lesssim  \left\|
|u|^{2\si}u\right\|_{L^{q'}(I;L^{r'})} 
\lesssim  \| u\|^{2\si}_{L^k(I;L^r)} \|u\|_{L^q(I;L^r)}\\
&\lesssim  \(\| v\|^{2\si}_{L^k(I;L^r)}+\|
w\|^{2\si}_{L^k(I;L^r)}\) \(\|v\|_{L^q(I;L^r)}+\|w\|_{L^q(I;L^r)}\), 
\end{aligned}
\end{equation}
for $I=[t_0,t]$. We have, from Lemma~\ref{lem:GN} and \eqref{eq:estimv}, 
\begin{equation*}
\|v(t)\|_{L^r}\lesssim \( \frac{1+\om_2}{\cosh (\om_1
t)}\)^{\delta(s)/n} = \( \frac{1+\om_2}{\cosh (\om_1
t)}\)^{\si/(2\si+2)}\, ,
\end{equation*}
therefore
\begin{equation}\label{eq:apparom}
\begin{aligned}
\|v\|_{L^k(I;L^r)}^k&\lesssim \int_{t_0}^t\( \frac{1+\om_2}{\cosh (\om_1
\tau)}\)^{2\si^2/(2-(n-2)\si)}d\tau \\
&\lesssim \(1+\om_2\)^{2\si^2/(2-(n-2)\si)}
\int_{t_0}^t \exp \(-\frac{2\om_1 \si^2 \tau}{2-(n-2)\si}\) d\tau \\
&\lesssim \frac{\(1+\om_2\)^{2\si^2/(2-(n-2)\si)}}{\om_1}
\exp \(-\frac{2\om_1 \si^2 t_0}{2-(n-2)\si}\)\\
&\lesssim \frac{\(1+\om_2\)^{2\si^2/(2-(n-2)\si)}}{\om_1}
\exp \(-\frac{2\om_1\si^2 \delta}{(2-(n-2))(1+\om_2)}\)
\, ,
\end{aligned}
\end{equation}
where we used Proposition~\ref{prop:local}. 
So long as $w$ satisfies
\eqref{eq:estimv}, the $L^qL^r$ norm of $w$ on the right hand side of
\eqref{eq:reste} can be absorbed,  provided that \eqref{eq:domin} is
satisfied for some sufficiently large $\Lambda$. Note that in general,
we cannot do without the last term of \eqref{eq:domin}, since for
$n\ge 2$ and $\si\ge \frac{2}{n}$, 
$\frac{2\si^2}{2-(n-2)\si}$ may be smaller or larger than $1$.  

We therefore get an estimate
for $\|w\|_{L^q(I;L^r)} $. Using Strichartz estimates again, 
\begin{equation*}
\begin{aligned}
\|w\|_{L^\infty(I;L^2)} &\lesssim  \left\|
|u|^{2\si}u\right\|_{L^{q'}(I;L^{r'})} \\
&\lesssim  \| u\|^{2\si}_{L^k(I;L^r)} \|u\|_{L^q(I;L^r)}\\
&\lesssim  \(\| v\|^{2\si}_{L^k(I;L^r)}+\|
w\|^{2\si}_{L^k(I;L^r)}\) \(\|v\|_{L^q(I;L^r)}+\|w\|_{L^q(I;L^r)}\)\\
&\lesssim
\(\frac{\(1+\om_2\)^{2\si^2/(2-(n-2)\si)}}{\om_1}\)^{2\si/k} 
e^{-\gamma \frac{\om_1}{1+\om_2}},  
\end{aligned}
\end{equation*}
for some positive $\gamma$. Applying the operator
$A$, we get similar estimates. We conclude by a  continuity
argument that $w$ satisfies \eqref{eq:estimv} for all $t\ge t_0$
provided that $\Lambda$ is sufficiently large in \eqref{eq:domin}, and
$u$ is global, with:
\begin{equation*}
Au\in L^\infty(\R;L^2)\cap L^q(\R;L^r)\, ,\quad \forall
A\in\{Id,J,H\}\, . 
\end{equation*}
Scattering follows
easily. Indeed, we have from Duhamel's principle:
\begin{equation*}
\left\| U_V(-\tau)u(\tau) -
U_V(-t)u(t)\right\|_{\Sigma} = \left\| \lambda \int_t^\tau
U_V(t-s) |u|^{2\si}u(s)ds \right\|_{\Sigma}\, .
\end{equation*}
From Strichartz and H\"older inequalities,  we have:
\begin{align*}
\left\| \int_t^\tau
U_V(t-s) |u|^{2\si}u(s)ds \right\|_{\Sigma}&\lesssim
\|u\|_{L^k([t,\tau];L^r)}^{2\si}\times\\
\times\sum_{A\in\{Id,
J,H\} } &
\left\|A u\right\|_{L^q([t,\tau];L^r)}\to 0 \quad \text{as }t,\tau\to
+\infty\, ,  
\end{align*}
since the $L^kL^r$ norm goes to $0$ exponentially, and the $L^qL^r$
norms are bounded. 

\section{Discussion}
\label{sec:discussion}

For the sake of simplicity, we discuss our results in the case
\begin{equation}\label{eq:Vhorse}
V(x) = \frac{1}{2}\(-\om_1^2 x_1^2 + \om_2^2 x_2^2\)\, .
\end{equation}
We proved global existence and scattering when $\si\ge \frac{2}{n}$, and
either $\|u_0\|_\Sigma$ is small, or \eqref{eq:domin} is satisfied. 
One expects that global existence (and scattering) should hold in
$\Sigma$ when the nonlinearity is defocusing, that is when $\lambda >0$. In
the absence of \emph{a priori} estimates ($E_V$ is not signed), we
cannot prove it. When $\si<\frac{2}{n}$ and $u_0\in \Sigma$, we proved
global existence, but not scattering. Of course, it holds when
$\|u_0\|_\Sigma$ is small, or \eqref{eq:domin} is satisfied.
What happens in the general case is not clear. On the other hand, proving
the existence of wave operators with asymptotic states in $\Sigma$ for
any $\si>0$ and any $\om_1>0$, $\om_2\ge 0$ is straightforward, thanks
to the exponential decay (see \cite{CaSIMA} for a proof which can be
easily adapted). We have more generally:
\begin{prop}
 Let $\lambda \in \R$, and $\si >0$, with $\si <\frac{2}{n-2}$
 if $n\geq 3$.    
Suppose that $V$ is of the form \eqref{eq:V}.
Then for every $u_- \in \Sigma$, there exist
 $T\in\R$ (finite) and a unique 
$u\in C(]-\infty,T];\Sigma)\cap L^q\( ]-\infty,T] ;
L^{2\si+2}(\R^n)\)$,
where $q=\frac{4\si+4}{n\si}$,  solution to
 \eqref{eq:nlsp} such that
$$\left\| U_V(-t)u(t)-u_-\right\|_\Sigma \Tend
t {-\infty} 0\, .$$
The same holds for positive times and $u_+\in \Sigma$. 
\end{prop}

The condition \eqref{eq:domin} is nonlinear with respect to $\om_2$,
which may seem surprising. On the other hand, notice that the constant
in front of $(1+\om_2)\ln (1+\om_2)$ is universal, unlike
$\Lambda$. This nonlinear factor vanishes when $\om_2=0$; it arises in
the estimate \eqref{eq:apparom}, and seems to be unavoidable. It is
due to the concentrations which are caused by the harmonic
potential (it is a linear phenomenon), and at time $t
=\frac{\pi}{2\om_2}$, this effect is not yet
balanced by the repulsive harmonic potential in the first
direction. \\

We do not know a criterion for global existence
(or, equivalently, finite time ``blow up'') other than boundedness in
$\Sigma$. One would expect the
unboundedness of $\|\nabla_x u(t)\|_{L^2}$ to be the only obstruction
to global existence. It turns out that this is true in the case of the
\emph{isotropic} repulsive harmonic potential studied in
\cite{CaSIMA}. Proving this was not straightforward though, and used a
particular evolution law, which does not seem to be available when the
potential is not isotropic. \\

When there is no potential, or when the potential is non-negative
(slightly more general potential are allowed, see \cite{CazCourant}),
a sufficient condition for finite time blow-up is provided by the
virial theorem, following the ideas of \cite{Glassey}. It relies in
the existence of a relatively simple evolution law for 
$\|xu(t)\|_{L^2}^2$. Sufficient conditions for finite time blow-up in
the case of the isotropic harmonic potential, or the isotropic
repulsive harmonic potential were obtained in a similar way in
\cite{CaAHP,CaSIMA}, using ordinary differential equations
techniques. In the case of the potential \eqref{eq:Vhorse}, it seems
that no such simplification is here to help. 
\bigskip

\noindent {\bf Acknowledgments}. The author is grateful to Jorge
Drumond Silva
for stimulating discussions about
Section~\ref{sec:strichartz}, and to the referee for pointing out
several imprecisions. Support by the European network HYKE, 
funded  by the EC as contract  HPRN-CT-2002-00282, is acknowledged.

\end{document}